
\input amstex
\documentstyle{amsppt}

\magnification=\magstep 1

\def\hak {^{\vee}}

\def \sgn {\operatorname{sgn}}

\document
\topmatter
\title{RESTRICTING SCHUBERT CLASSES}
\endtitle
\author{Piotr Pragacz}
\endauthor
\date July 31, 1999
\enddate
\address Mathematical Institute of Polish Academy of Sciences, 
Chopina 12, PL 87100 Toru\'n, Poland
\endaddress
\thanks{This note was supported by KBN grant No. 2P03A 05112.}
\endthanks
\endtopmatter

\noindent
{\bf Introduction}
\smallskip

The goal of the present note is twofold. Firstly, we correct some points  
in the paper (quoted in the following as [P]):
\block
P. Pragacz:  ``A generalization of the Macdonald-You formula",
Journal of Algebra {\bf 204}, 573--587 (1998).
\endblock
Secondly, by comparing the formula appearing in the title of [P] with
results of Stembridge [St] and some other 
combinatorial results, we deduce some new identities
in Propositions 2--6. They concern: restrictions of Schubert classes
to the cohomology of Lagrangian Grassmannians as well as relations
between $Q$-functions, Stembridge's coefficients, and various 
``hook numbers". Also, we provide some examples illustrating [P] and 
the formulas given in the present note. 
\smallskip

In this note, any unexplained notation or quotation stems 
from [P]. However, in order to  make the notation
maximally compatible with that used in [St] (which is our prinicipal
reference here), we label strict partitions by $\lambda$,
and $\mu$ usually denotes an ordinary partition, contrariwise to [P].

\bigskip

\noindent
{\bf 1. Erratum to [P]}
\smallskip

Due to some bugs in the computer system SCHUR [Sch], 
[P, Example 3(b)] was
miscalculated: the quadratic expression in $Q$-functions displayed there, 
written as a $\Bbb Z$-linear combination of $Q$-functions, 
contains {\it no} negative summands.
Consequently the sentence on p.585, lines 6--7 from the bottom, is to
be withdrawn from [P].
(These corrections {\it do not} affect other results of [P], in particular
the main formulas.)

\bigskip

\noindent
{\bf 2. Nonnegativity of the restriction coefficients}
\smallskip

In fact, if 
$$
i^*(\sigma_{\mu})=\sum_{\lambda} c_{\lambda \mu} 
\sigma'_{\lambda},
\tag 1
$$
with $c_{\lambda \mu}\in \Bbb Z$, then all the coefficients $c_{\lambda \mu}$ 
are nonnegative. 
Perhaps the easiest way to see this, is the following.
Let for $a\in H^*(G;\Bbb Z)$,  
$\int_G a$ stand for
the degree of the top codimensional component of $a$, and define similarly
$\int_{G'} b$ \  for $b\in H^*(G';\Bbb Z)$.
Given a strict partition $\lambda\subset (n,{n-1},\ldots,1)$, we denote by
$\lambda\hak$ the strict partition whose parts complement those
of $\lambda$ in $\{1,\ldots,n\}$. We record the following property [P2]:

\proclaim{Lemma 1 (Duality)} \ The basis $\{\sigma'_{\lambda}\}$ of the group 
$H^{2p}(G';\Bbb Z)$ and the basis $\{\sigma'_{\lambda\hak}\}$ of the group
$H^{n(n+1)-2p}(G';\Bbb Z)$ are dual under the pairing \ $(a,b) \mapsto
\int_{G'} a\cdot b$ of Poincar\'e duality. 
\endproclaim

Now, if $i^*(\sigma_{\mu})=\sum_{\lambda} c_{\lambda \mu} \sigma'_{\lambda}$,
with $c_{\lambda \mu}\in \Bbb Z$, then it follows from the duality 
property that 
$$
c_{\lambda \mu}=\int_{G'} i^*(\sigma_{\mu})\cdot \sigma'_{\lambda\hak}.
\tag 2
$$
Using the projection formula for $i$, this is rewritten as
$$
c_{\lambda \mu}=\int_G \sigma_{\mu}\cdot i_*(\sigma'_{\lambda\hak}).
\tag 3
$$
Regard $G$ as a homogeneous space $GL(V)/P$, where $P$ is a suitable
parabolic subgroup of $GL(V)$. Let $\Omega\subset G$ be a Schubert
variety representing $\sigma_{\mu}$ and let $\Omega'\subset G'
\subset G$ be a Schubert variety representing $\sigma'_{\lambda\hak}$.
Using e.g.
Kleiman's theorem on a general translate [K], we can replace $\Omega$
by a translate by an element $g\in GL(V)$ such that $g\cdot \Omega$ 
and $\Omega'$ meet properly, and this intersection is represented
as a nonnegative zero-cycle. This shows that $c_{\lambda \mu}\ge 0$.

\smallskip

A similar property holds in the following more general setting. 
Let now \  $G\supset P\supset B$ be a semisimple linear algebraic group, 
a parabolic subgroup, and a Borel subgroup.
In a generalized flag variety $G/P$, one has Schubert varieties
$\overline{BwP/P}$ and their Schubert classes in $H^*(G/P;\Bbb Z)$
indexed by a corresponding subset of the Weyl group. These Schubert
classes enjoy a similar duality property. 
In an analogous way, using a general translate argument, one shows
that the fundamental class of any subscheme of 
$G/P$ is a $\Bbb Z$-linear combination of the Schubert classes in 
$H^*(G/P;\Bbb Z)$ with nonnegative coefficients. Combining this with
a well-known fact about pulling back the class of a Cohen-Macaulay
subscheme (see, e.g., Lemma on p.108 in [F-P]), we get the following 
result (also implying the nonnegativity of the above $c_{\lambda \mu}$):

\proclaim{Proposition 1} Let $f: G/P \to Y$ be morphism to a nonsingular
variety $Y$. 
Let $Z$ be a pure-dimensional closed Cohen-Macaulay subscheme of $Y$. 
Then $f^*([Z])$ is a $\Bbb Z$-linear combination of the Schubert classes 
in $H^*(G/P;\Bbb Z)$ with nonnegative coefficients.
\endproclaim 

\bigskip

\noindent
{\bf 3. Stembridge's coefficients}
\smallskip

We recall (see the discussion after [P, Proposition 8]) that the
coefficients appearing in (1) and those appearing in:
$$
\eta(s_{\mu})=\sum_{\lambda} g_{\lambda \mu} Q_{\lambda}
\tag 4
$$
satisfy \ $c_{\lambda \mu}=g_{\lambda \mu}$. Here, we take suficiently 
large Grassmannians $i: G'\hookrightarrow G$. To be more precise, this 
means that given $\mu$, we take $n\ge |\mu|$ so that any strict partition
$\lambda$ with $|\lambda|=|\mu|$ is contained in $(n,n-1,\ldots,1)$.
Consequently, all the coefficients $g_{\lambda \mu}$ are nonnegative.
But this result, together with a combinatorial interpretation of
the $g_{\lambda \mu}$'s, was already established by Stembridge in [St].
\footnote {The fact that this result was already established by Stembridge,
has been learned by the author only in June 1999.}
Indeed, the last displayed (unnumbered) equality before [St, Theorem 9.3]: 
$$
`` \ S_{\mu}=\sum_{\lambda \in DP_n} g_{\lambda \mu} Q_{\lambda} \ "
\tag 5
$$ 
is identical with (4) because  $S_{\mu}$ in the notation of [St] (and [M])
is equal to $\eta(s_{\mu})$ in our notation. \footnote {Note that 
the map denoted in [P] and here by $\eta$, is denoted by $\varphi$ in [M].} 
In [St], (5) is a consequence of the equality 
$$
P_{\lambda}=\sum_{|\mu|=|\lambda|} g_{\lambda \mu} s_{\mu},
\tag 6
$$
where $P_{\lambda}=2^{-l(\lambda)} Q_{\lambda}$,
and comparison of the canonical scalar products on the ring of all 
symmetric functions
with that on the ring of $Q$-functions. To the nonnegativity of 
$g_{\lambda \mu}$
is given, in loc. cit., several interpretations in representation theory,
some of which go back to Morris and Stanley. 

(Observe that (4) and (6) yield the following expression for 
$\eta(P_{\lambda})$:
$$
\eta(P_{\lambda})=\sum_{|\mu|=|\lambda|} g_{\lambda \mu} \eta(s_{\mu})
=\sum_{|\mu|=|\lambda|} \sum_{|\nu|=|\lambda|} g_{\lambda \mu} g_{\nu \mu}
Q_{\nu}.)
\tag 7
$$

Stembridge [St] also established a combinatorial interpretation
of the numbers $f_{\mu \nu}^{\lambda}$ appearing as coefficients in
the expansion:
$$
P_{\mu}P_{\nu}=\sum_{\lambda} f_{\mu \nu}^{\lambda} P_{\lambda},
\tag 8
$$
where $\mu$, $\nu$, and $\lambda$ denote now strict partitions. It will be 
convenient to set 
$$
e_{\mu \nu}^{\lambda}:= 2^{l(\mu)+l(\nu)-l(\lambda)} f_{\mu \nu}^{\lambda}.
\tag 9
$$
There exists a geometric analogue of (8): in the cohomology
ring $H^*(G'; \Bbb Z)$ of a sufficiently large Lagrangian Grassmannian,
$$
\sigma'_{\mu}\cdot \sigma'_{\nu}= \sum_{\lambda} 
e_{\mu \nu}^{\lambda} \sigma'_{\lambda}.
\tag 10
$$
(See [P2, Sect.6].)
\smallskip

Stembridge's combinatorial description of the above $f_{\mu \nu}^{\lambda}$ 
and $g_{\lambda \mu}$ can be summarized by the following:

\proclaim{Theorem [St]} \ (i) The coefficient $f_{\mu \nu}^{\lambda}$
is equal to the number of marked shifted tableaux $T$ of shape $\lambda/\mu$
and weight (or content) $\nu$ such that:

\noindent
(a) The word w(T) associated with $T$ ([St, Sect.8] and [M, p.258]) 
has the lattice property in the sense of loc.cit.;

\noindent
(b) for each $k\ge 1$, the rightmost occurence of $k'$ in $w(T)$ precedes
the last occurence of $k$. 
\smallskip
\noindent
(ii) The coefficient $g_{\lambda \mu}$ is equal to the number of unshifted
marked tableaux $T$ of shape $\mu$ and weight $\lambda$ satisfying (a)
and (b) above.
\endproclaim
For all unexplained here combinatorial notions, we refer the reader
to [St, Sect.6 and 8], [P2, Sect.4], and to [M, III.8 pp.255--259]. 
We make no attempt to make a complete survey here. 
Some examples of the coefficients 
$g_{\lambda \mu}$ will be given below.

\smallskip

Summarizing the content of this section, we record:

\proclaim{Proposition 2} \ We have for a partition $\mu\subset (n^n)$
$$
i^*(\sigma_{\mu}) = \sum_{\lambda} g_{\lambda \mu} \ \sigma'_{\lambda} \,,
\tag 11
$$
where $\lambda$ runs over strict partitions contained in $(n,n-1,\ldots,1)$, 
and
$g_{\lambda \mu}$ is the Stembridge coefficient described in Theorem (ii).
\endproclaim

\bigskip

\noindent
{\bf 4. Quadratic relations between $Q$-functions}
\smallskip

We pass now to some applications of the generalized Macdonald-You formula
([L-L2], [P, Corollary 2]):
$$
2^n \eta(s_{\mu})=\sum Q_{(a_{i_1},\ldots,a_{i_k})}\cdot
Q_{A \# B \smallsetminus (a_{i_1},\ldots,a_{i_k})}.
\tag 12
$$
Recall that here, for $\mu=(\alpha_1,\ldots,\alpha_n|
\beta_1,\ldots,\beta_n)$ in Frobenius notation,
$$
A=(a_1,\ldots,a_n):=(\alpha_1+1,\ldots,\alpha_n+1)\,, \ \ \ \ \ 
B:=(\beta_1,\ldots,\beta_n)\,,
\tag 13
$$ 
and the sum is over all sequences $1\le i_1<\cdots < i_k\le n$ 
and $k=0,1,\ldots, n$.

\smallskip
 
Since $\eta(e_i)=\eta(h_i)$, where $h_i$ is the $i$th complete homogeneous
symmetric function, we have for a partition $\mu$ 
$$
\eta(s_{\mu^{\sim}})=\eta(s_{\mu}),
\tag 14
$$ 
where $\mu^{\sim}=
(\beta_1,\ldots,\beta_n|
\alpha_1,\ldots,\alpha_n)$ is the conjugate partition of $\mu$.
We set in addition 
$$
C=(c_1,\ldots,c_n):=(\beta_1+1,\ldots,\beta_n+1)\,, \ \ \ \ \
D:=(\alpha_1,\ldots,\alpha_n).
\tag 15
$$
Then (12) and (14) imply the following:

\proclaim{Proposition 3} \ We have  
$$
\sum Q_{(a_{i_1},\ldots,a_{i_k})}\cdot
Q_{A \# B \smallsetminus (a_{i_1},\ldots,a_{i_k})}
=\sum Q_{(c_{i_1},\ldots,c_{i_k})}\cdot
Q_{C \# D \smallsetminus (c_{i_1},\ldots,c_{i_k})},
\tag 16
$$
where the sums are over all sequences $1\le i_1<\cdots <i_k\le n$ 
and $k=0,1,\ldots, n$. 
\endproclaim

The relations (16), regarded from the side of $Q$-functions, seem to be rather 
nontrivial.
For instance, for $\mu=(5^331^3)=(432|621)$, so $A=(5,4,3)$, $B=(6,2,1)$, 
$C=(7,3,2)$, and $D=(4,3,2)$, we get the equation
$$\aligned
&Q_{654321} - Q_5\cdot Q_{64321} + Q_4\cdot Q_{65321} - Q_3\cdot Q_{65421} 
- Q_{54}\cdot Q_{6321} \\
&+ Q_{53}\cdot Q_{6421} - Q_{43}\cdot Q_{6521} + Q_{543}\cdot Q_{621} \\
&= Q_{32}\cdot Q_{7432} + Q_{732}\cdot Q_{432}. 
\endaligned
\tag 17
$$
Using a (hopefully) debugged version of SCHUR, (16) is expressed as the 
following $\Bbb Z$-linear combination of the $Q_{\lambda}$'s :
\setbox1=\vbox{\settabs5\columns{
\+$8Q_{11\ 64}$&$ + \ 8Q_{11\ 631}$&$ + \ 8Q_{11\ 541}$&$ + \ 8Q_{11\ 532}$&$ + \ 8Q_{10\ 74}$\cr
      \+$ + \ 8Q_{10\ 731}$&$ + \ 8Q_{10\ 65}$&$ + \ 24Q_{10\ 641}$&$ + \ 24Q_{10\ 632}$&$ + \ 24Q_{10\ 542}$\cr
      \+$ + \ 8Q_{10\ 5321}$&$ + \ 8Q_{975}$&$ + \ 16Q_{9741}$&$ + \ 16Q_{9732}$&$ + \ 16Q_{9651}$\cr
      \+$ + \ 48Q_{9642}$&$ + \ 16Q_{96321}$&$ + \ 16Q_{9543}$&$ + \ 16Q_{95421}$&$ + \ 8Q_{8751}$\cr
      \+$ + \ 24Q_{8742}$&$ + \ 8Q_{87321}$&$ + \ 24Q_{8652}$&$ + \ 24Q_{8643}$&$ + \ 24Q_{86421}$\cr
      \+$ + \ 8Q_{85431}$&$ + \ 8Q_{7653}$&$ + \ 8Q_{76521}$&$ + \ 8Q_{76431}$\cr}}
$$\box1$$
Consequently, taking sufficiently large Grassmannians $i: G'\hookrightarrow G$,
we have 
\smallskip

$i^*(\sigma_{5^331^3})=$
\setbox1=\vbox{\settabs5\columns{
\+$\sigma'_{11\ 64}$&$ + \ \sigma'_{11\ 631}$&$ + \ \sigma'_{11\ 541}$&$ + \ \sigma'_{11\ 532}$&$ + \ \sigma'_{10\ 74}$\cr
      \+$ + \ \sigma'_{10\ 731}$&$ + \ \sigma'_{10\ 65}$&$ + \ 3\sigma'_{10\ 641}$&$ + \ 3\sigma'_{10\ 632}$&$ + \ 3\sigma'_{10\ 542}$\cr
      \+$ + \ \sigma'_{10\ 5321}$&$ + \ \sigma'_{975}$&$ + \ 2\sigma'_{9741}$&$ + \ 2\sigma'_{9732}$&$ + \ 2\sigma'_{9651}$\cr
      \+$ + \ 6\sigma'_{9642}$&$ + \ 2\sigma'_{96321}$&$ + \ 2\sigma'_{9543}$&$ + \ 2\sigma'_{95421}$&$ + \ \sigma'_{8751}$\cr
      \+$ + \ 3\sigma'_{8742}$&$ + \ \sigma'_{87321}$&$ + \ 3\sigma'_{8652}$&$ + \ 3\sigma'_{8643}$&$ + \ 3\sigma'_{86421}$\cr
      \+$ + \ \sigma'_{85431}$&$ + \ \sigma'_{7653}$&$ + \ \sigma'_{76521}$&$ + \ \sigma'_{76431}$\cr}}
$$\box1$$
So e.g. we have:  \ $g_{(5^331^3)\ (11\ 64)}=1, 
\ g_{(5^331^3)\ (10\ 641)}=3, 
\ g_{(5^331^3)\ (9741)}=2$, and \ $g_{(5^331^3)\ (9642)}=6$. 

\bigskip

\noindent
{\bf 5. Linear relations between Stembridge's coefficients}
\smallskip

Combining (4), (12), and (16), we have in the above notation, 
associated with a fixed $\mu$
$$\aligned
&\sum Q_{(a_{i_1},\ldots,a_{i_k})}\cdot
Q_{A \# B \smallsetminus (a_{i_1},\ldots,a_{i_k})} \\
&=\sum Q_{(c_{i_1},\ldots,c{i_k})}\cdot
Q_{C \# D \smallsetminus (c_{i_1},\ldots,c_{i_k})} 
= 2^n \sum_{\lambda} g_{\lambda \mu} Q_{\lambda},
\endaligned
\tag 18
$$
where the first two sums are over all sequences $1\le i_1<\cdots <i_k\le n$ 
and $k=0,1,\ldots, n$.

The equalities (18) imply linear relations between the 
$e_{\mu \nu}^{\lambda}$'s and $g_{\lambda \mu}$'s.
Given a sequence of different positive integers $K=(k_1,\ldots,k_l)$,
there is a permutation $w=w_K \in S_l$ such that $k_{w(1)} >
\cdots > k_{w(l)} > 0$. Denote this last-mentioned strict partition
by $<K>$. Then given strict partitions $\mu$, $\lambda$
and  a sequence $K$ as above, we set 
$$
e_{\mu \ K}^{\lambda}:= \sgn (w_K) \ e_{\mu \ <K>}^{\lambda}.
\tag 19
$$
From (18) and (16) we get the following result:

\proclaim{Proposition 4} \ For a fixed partition $\mu$ and strict partition 
$\lambda$ with $|\mu|=|\lambda|$, we have in the above notation associated
with $\mu$
$$\aligned
2^n g_{\lambda \mu} &= \sum e^{\lambda}_{(a_{i_1},\ldots, a_{i_k}), \ 
A\#B \smallsetminus (a_{i_1},\ldots,a_{i_k})} \\
&= \sum e^{\lambda}_{(c_{i_1},\ldots, c_{i_k}), \ 
C\#D \smallsetminus (c_{i_1},\ldots,c_{i_k})}
\endaligned
\tag 20
$$
where the sums are over all sequences $1\le i_1<\cdots <i_k\le n$ 
for which $A\#B \smallsetminus (a_{i_1},\ldots,a_{i_k})$ (resp. 
$C\#D \smallsetminus (c_{i_1},\ldots,c_{i_k})$) is a sequence of
different integers, and $k=0,1,\ldots, n$.
\endproclaim
For instance, for any strict partition $\lambda$ with $|\lambda|=21$,
and for $\mu=(5^331^3)=(432|621)$, we get the equations:
$$\aligned
2^3 g_{\lambda \ (5^331^3)} \ = \ & e^{\lambda}_{(654321) \  (\emptyset)} 
- e^{\lambda}_{(5) \ (64321)} + e^{\lambda}_{(4) \ (65321)} 
- e^{\lambda}_{(3) \ (65421)} - e^{\lambda}_{(54) \ (6321)} \\ 
&+ e^{\lambda}_{(53) \ (6421)} - e^{\lambda}_{(43) \ (6521)} 
+ e^{\lambda}_{(543) \ (621)} \\ 
\ = \ &e^{\lambda}_{(32) \ (7432)} + e^{\lambda}_{(732) \ (432)}. 
\endaligned
$$

\bigskip

\noindent
{\bf 6. The class $i_*(\sigma'_{\lambda})$ as a $\Bbb Z$-linear combination
of the $\sigma_{\mu}$'s}
\smallskip

By reasoning similarly as in $\S$4, one shows that
for any proper morphism $f: X \to G/P$ from a scheme $X$ to a generalized
flag variety $G/P$, and for any irreducible subscheme $Y\subset X$,
$f_*([Y])$ is a $\Bbb Z$-linear combination of Schubert classes
in $H^*(G/P; \Bbb Z)$ with nonnegative coefficients.  

\smallskip

The next proposition will give $i_*(\sigma'_{\lambda})$ as an 
explicit $\Bbb Z$-linear
combination of the $\sigma_{\mu}$'s. Given a partition $\mu\subset (n^n)$,
we set $\mu^{\star}:=(n-{\mu_n},\ldots,n-{\mu_1})$. The following duality
property is a well-known result of Schubert calculus [F]:

\proclaim{Lemma 2} \ 
The basis $\{\sigma_{\mu}\}$ of the group 
$H^{2p}(G;\Bbb Z)$ and the basis $\{\sigma_{\mu^{\star}}\}$ of the group
$H^{2(n^2-p)}(G;\Bbb Z)$ are dual under the pairing \ $(a,b) \mapsto
\int_{G} a\cdot b$ of Poincar\'e duality. 
\endproclaim

We now state:

\proclaim{Proposition 5} \ For a fixed strict partition $\lambda 
\subset (n,n-1,\ldots, 1)$, we have
$$
i_*(\sigma'_{\lambda}) = \sum_{|\mu|=|\lambda|+n(n-1)/2} 
g_{\lambda\hak, \mu^{\star}} \  
\sigma_{\mu}\,,
\tag 21
$$
where $\mu$ runs over partitions contained in $(n^n)$ and 
$g_{\lambda\hak, \mu^{\star}}$
is the Stembridge coefficient described in Theorem (ii).
\endproclaim
Indeed, if $i_*(\sigma'_{\lambda})=\sum_{\mu} m_{\lambda \mu} 
\sigma_{\mu}$, with $m_{\lambda \mu}\in \Bbb Z$ (so that
$|\mu|=|\lambda|+n(n-1)/2$ \ ), then  
it follows from Lemma 2 that 
$$
m_{\lambda \mu}=\int_{G} (i_* \sigma'_{\lambda}) \cdot \sigma_{\mu^{\star}}.
\tag 22
$$
Using the projection formula for $i$, this is rewritten as
$$
m_{\lambda \mu}=\int_{G'} \sigma'_{\lambda}\cdot i^*(\sigma_{\mu^{\star}}).
\tag 23
$$
In turn, using the description of $i^*(\sigma_{\mu^{\star}})$ from 
Proposition 2,
(23) is rewritten as
$$
m_{\lambda \mu} = \int_{G'} \sigma'_{\lambda}\cdot (\sum_{\nu} 
g_{\nu \mu^{\star}} \ \sigma'_{\nu})
= \int_{G'}\sum_{\tau} \sum_{\nu} e_{\lambda \nu}^{\tau} \ g_{\nu \mu^{\star}}
 \ \sigma'_{\tau}=g_{\lambda\hak, \mu^{\star}}
\tag 24
$$
because only $\tau=(n,n-1,\ldots,1)$ and $\nu=\lambda\hak$ give a nonzero
contribution (note that for such $\tau$ and $\nu$, we have
$e_{\lambda \nu}^{\tau}=1$). 

\bigskip
\noindent
{\bf 7. Relations between the degrees of the ordinary and projective
representations of the symmetric groups}
\smallskip

For a partition $\mu$, we set
$$
\overline f^{\mu} := \prod_{x\in {\mu}} {1\over {h(x)}},
\tag 25
$$
where $h(x)$ is the hook-length of $\mu$ at $x=(i,j)$ defined by
$h(x)=h(i,j)=\mu_i + \mu^{\sim}_j -i -j +1$. If $|\mu|=m$ then
$f^{\mu}:= m! \ {\overline f^{\mu}}$ is the degree of the irreducible
representation of $S_m$ corresponding to $\mu$. Equivalently, $f^{\mu}$
is the number of standard tableaux of shape $\mu$, obtained by labeling
the squares of the diagram of $\mu$ with the numbers $1,2,\ldots,m$. 
We refer to [F] for a detailed discussion of these facts.

\smallskip

For a strict partition $\lambda$, we set
$$
\overline g^{\lambda} := \prod_{x\in S({\lambda})} {1\over {h(x)}},
\tag 26
$$
where $S(\lambda)$ is the shifted diagram associated with $\lambda$
[M, p.255], and for each square $x\in S(\lambda)$ the hook-length $h(x)$
is defined to be the hook-length at $x$ in the ``double diagram"
$(\lambda_1,\lambda_2,\ldots |\lambda_1-1,\lambda_2-1,\ldots)$, 
containing $S(\lambda)$.
If $|\lambda|=m$, $g^{\lambda}:= m! \ {\overline g^{\lambda}}$ is the number
of shifted standard tableaux of shape $S(\lambda)$, obtained by labeling
the squares of $S(\lambda)$ with the numbers $1,2,\ldots,m$ with strict
increase along each row and down each column. The numbers $g^{\lambda}$
also admit an interpretation as the degree of suitable projective
representations of $S_m$. We refer to [H-H] for a detailed discussion
of these results.

One has the following formulas, in terms of parts, for 
$\overline f^{\mu}$ [M, I.1 Example 1]
and $\overline g^{\lambda}$ [M, III.8 Example 12]:
$$
\overline f^{\mu} = {{ \prod_{i<j} (\mu_i - \mu_j - i + j) } \over
{\prod_{i\ge 1}(\mu_i+n-i)!}} \,,
\tag 27
$$
$$
\overline g^{\lambda} = {1 \over {\prod_{i\ge 1}{\lambda_i} !}} \ 
{\prod_{i<j} {{\lambda_i - \lambda_j} \over {\lambda_i + \lambda_j}}} .
\tag 28
$$

We now record:

\proclaim{Lemma 3} \ (i) Under the specialization $e_i:= {1 \over {i!}}$ , 
$s_{\mu}$ becomes \ $\overline f^{\mu}$.
\smallskip
\noindent
(ii) Under the specialization $Q_i:= {1 \over {i!}}$ , $Q_{\lambda}$ becomes
 \ $\overline g^{\lambda}$.
\endproclaim
\noindent
(For assertion (i), see [M, I.3 Example 5]. Assertion (ii) stems from
[DC-P, Proposition 6].) 

\smallskip

Given a partition $\mu$, we want to apply formulas (12) and (16), so
we adopt the notation of $\S$5. Also, we follow the notation of
$\S$6 associated with a sequence $K$. For such a sequence, we set
$$
\overline g^K:= \sgn (w_K) \ \overline g^{<K>}.
\tag 29
$$
 
From Lemma 3, (12), and (16),
we get

\proclaim{Proposition 6} \ For a fixed partition $\mu$, we have 
$$
\aligned
2^n \overline f^{\mu} &= \sum \overline g^{(a_{i_1},\ldots, a_{i_k})} \ 
\overline g^{A\#B \smallsetminus (a_{i_1},\ldots,a_{i_k})} \\
&=\sum \overline g^{(c_{i_1},\ldots, c_{i_k})} \ 
\overline g^{C\#D \smallsetminus (c_{i_1},\ldots,c_{i_k})}\,,
\endaligned
\tag 30
$$
where the sums are over all sequences $1\le i_1<\cdots <i_k\le n$ 
for which $A\#B \smallsetminus (a_{i_1},\ldots,a_{i_k})$ (resp. 
$C\#D \smallsetminus (c_{i_1},\ldots,c_{i_k})$) is a sequence of
different integers, and $k=0,1,\ldots, n$.
\endproclaim

For instance, for $\mu=(5^331^3)=(432|621)$, we get the equations:
$$\aligned
2^3 \overline f^{(5^331^3)} \ = \ & \overline g^{(654321)} - 
\overline g^{(5)} \ \overline g^{(64321)} + \overline g^{(4)} \ 
\overline g^{(65321)} - \overline g^{(3)} \ 
\overline g^{(65421)} - \overline g^{(54)} \ \overline g^{(6321)} \\ 
&+ \overline g^{(53)} \ \overline g^{(6421)}  
-\overline g^{(43)} \ \overline g^{(6521)} + \overline g^{(543)} \ 
\overline g^{(621)} \\ 
 \ = \ &\overline g^{(32)} \ \overline g^{(7432)} + \overline g^{(732)} \ 
\overline g^{(432)}. 
\endaligned
$$

\vskip 10pt
\noindent
{\eightrm ACKNOWLEDGEMENTS. \ The author would like to thank Brian Wybourne
for debugging, on the author's request, 
the computer program SCHUR in June 1999. 
On the occasion of this debugging several examples were computed, which
have helped the author to realize the link of [P]
with Stembridge's results on shifted tableaux.

We also thank Anders Thorup for his advice in preparing this note.}

\enddocument
\bye